\documentclass{ieeetran}
\usepackage{cite}
\usepackage{amsmath,amssymb,amsfonts}
\usepackage{algorithmic}
\usepackage{graphicx}
\usepackage{textcomp}
\usepackage{cite}
\usepackage{amsmath,amssymb,amsfonts}
\usepackage{algorithmic}
\usepackage{graphicx}
\usepackage{textcomp}

\usepackage{epsfig}
\usepackage{latexsym}
\usepackage{pdfsync}

\newcommand{\bc}{\begin{center}}
\newcommand{\ec}{\end{center}}

\newcommand{\eq}{\begin{equation}\begin{array}{rllllllllllllllllllllllllllllllll}}
\newcommand{\ee}{\end{array}\end{equation}}
\newcommand{\bmt}{\left[ \begin{array}{ccccccccc}}
\newcommand{\emt}{\end{array}\right]}
\newcommand{\bea}{\begin{eqnarray}}
\newcommand{\eea}{\end{eqnarray}}
\newcommand{\bean}{\begin{eqnarray*}}
\newcommand{\eean}{\end{eqnarray*}}

\begin{document}
\title{Optimal Boundary Control of a Nonlinear Reaction Diffusion Equation via Completing \\ the Square and Al'brekht's Method}
\author{Arthur J. Krener, \IEEEmembership{Life Fellow}},
\thanks{This work was supported by AFOSRunder FA9550-20-1-0318.}
\thanks{A. J. Krener is with the Department of Applied Mathematics, Naval Postgraduate School,
Monterey, CA 93940, USA, 
        {\tt\small ajkrener@nps.edu}}
\maketitle
\begin{abstract}
 The two contributions of this paper are as follows.  The first is the solution 
of an  infinite dimensional,  boundary controlled Linear Quadratic Regulator by the 
simple and constructive method of completing the square.  The second contribution
is the extension of Al'brekht's method to the optimal stabilization of a boundary
controlled, nonlinear Reaction Diffusion system.
\end{abstract}


\section{Introduction}
In 1961 Al'brekht \cite{Al61} showed how one could compute degree by degree the Taylor polynomials
of the optimal cost and optimal feedback of a smooth, nonlinear, infinite horizon, finite dimensional  optimal control problem
provided the linear part of the dynamics and the quadratic part of the running cost satisfied the standard
Linear Quadratic Regulator (LQR) conditions.

 Recently \cite{Kr20} we showed  how Al'brekht's method could be generalized to infinite dimensional problems
 with distributed control.  In this paper we show  how Al'brekht's method can be generalized to infinite dimensional problems
 with boundary  control.     In the next section we present and explicitly solve an LQR for the boundary control of a heated rod.  
 We do this in a novel way, by completing the square in infinite dimensions.  In Section Three we analyze the closed loop linear dynamics.
 In Section Four
we show how Al'brekht's method can be used to stabilize a nonlinear reaction diffusion equation using boundary control.

We are not the first to use Al'brekht's method on infinite dimensional systems, see the works of   Kunisch and coauthors \cite{BKP18}, \cite{BKP19}, \cite{KP19}.  Krstic, Vazquez and coauthors have had great success stabilizing infinite dimensional systems through boundary control where  the nonlinearites are expressed  by Volterra  integral operators of increasing degrees using backstepping techniques, \cite{KS08}, \cite{VK08}.   In our  extension of Al'brekht, we assume that the nonlinearites are given by Fredholm integral expressions of increasing degrees.

\section{ Boundary Control of a Heated   Rod}

We consider a modification of Example 3.5.5 of Curtain and Zwart \cite{CZ95}.
We have a rod of length one insulated at one end and heated/cooled
at the other.  The goal is to control the temperature to a constant set point which we conveniently take to be zero.

Let $0\le x\le 1$ be distance along the rod, $z(x,t)$ be the temperature 
of the rod at $x,t$ and $z^0(x)$ be the initial temperature distribution of the 
rod at $t=0$.  The goal is to stabilize the temperature to $z=0$ as $t\to \infty$
using boundary control at $x=1$.  

The rod is modeled by these equations
\bea
\frac{\partial z}{\partial t}(x,t)&=&\frac{\partial^2 z}{\partial x^2}(x,t) \label{old}\\
z(x,0)&=& z^0(x) \label{ic}\\
\frac{\partial z}{\partial x}(0,t)&=&0 \label{BC0}\\
\frac{\partial z}{\partial x}(1,t)&=& \beta (u(t)-z(1,t)) \label{BC1}
\eea
for some positive constant $\beta $ where the control is $u(t)$, the temperature applied to the end of the rod.

First  we study the open loop system where $u(t)=0$ for all $t\ge 0$.   We consider the closed linear operator
\bean
Ah(x) &=&\frac{d^2 h }{d x^2} (x)
\eean
 where its  domain is space of all $h\in {\bf L}_2[0,1]$ such that $h$ and $\frac{dh}{dx} $ are absolutely continuous
 and $\frac{dh}{dx}(0)=0$, $\frac{dh}{dx}(1)=-  h(1)$.
Because of the Neumann boundary condition at $x=0$, 
the eigenvectors are of form
\bean
\phi(x) &=&c \cos \nu x
\eean
for some constants $ \nu$ and $c\ne 0$.
The Robin boundary condition at $x=1$ implies that 
  $\nu $ is a root of the equation
  \bea \label{nueq1}
  \nu \sin \nu&=&\beta  \cos \nu
  \eea
  or equivalently
\bea  \label{nueq}
\nu&=&\beta  \cot \nu
\eea

There is one root, $\nu_n$, of this equation in each open interval $(n\pi, (n+1/2)\pi)$
for $n=0,1, \ldots$. The  $n^{th}$ root $\nu_n\to n\pi$ as $n\to \infty$.   As $ \beta  \to 0 $
the $n^{th}$ root $\nu_n\to n\pi$ and as $ \beta  \to \infty$ the $n^{th}$ root $\nu_n\to (n+1/2)\pi$. If $\beta =0$, $\nu_n= n\pi$ and we have an uncontrolled rod with no heat flux at either end.  
As $\beta  \to \infty$ the $n^{th}$ root $\nu_n\to (n+1/2) \pi$.
The corresponding eigenvalues are $\lambda_n=-\nu_n^2$.   

 If $  \beta=1$,  the first five roots are  $\nu_0\approx0.8603$, $\nu_1\approx3.4256$,  $\nu_2\approx6.4373$, $\nu_3\approx9.5293$ and $\nu_4\approx12.6453$. 
 So the five least stable eigenvalues are $\lambda_0=-\nu_0^2\approx-0.7402$, $\lambda_1=-\nu_1^2\approx-11.7349$, $\lambda_2=-\nu_2^2\approx-41.4388$,
 $\lambda_3=-\nu_3^2\approx-90.8082$ and $\lambda_4=-\nu_4^2\approx-159.9033$.  Notice that as $n\to \infty$, $\nu_n $ is monotonically decreasing to $n\pi$ and $\lambda_n $ is monotonically increasing to $-n^2\pi^2$.

   Because the Laplacian
is self adjoint with respect to these boundary conditions, the eigenfunctions
are orthogonal.  We  normalize them
\bea
\phi_n(x)&=&c_n\cos \nu_n x \label{ef}
\eea 
where
\bea \label{cn}
c_n&=&\mbox{sign}(\cos \nu_n)\ \sqrt{{4\nu_n \over 2\nu_n +\sin 2 \nu_n}}\label{cn}
\eea 
to get an orthonormal family satisfying $\phi_n(1)>0$.  Because $\nu_n\in (n\pi, (n+1/2)\pi)$
it follows that $\mbox{sign}\ c_n=\mbox{sign}(\cos \nu_n)=(-1)^n$ and $2\nu_n\in (2n\pi, (2n+1)\pi)$  so $\sin 2 \nu_n$
is positive.
This implies that  $|c_n|\le \sqrt{2}$ so
\bea\label{phn1}
\phi_n(1)&=&c_n\cos \nu_n \ \le \sqrt{2} \ |\cos \nu_n| \ \le \ \sqrt{2}
\eea 
Since $\nu_n-n\pi>0$ is monotonically decreasing to $0$ as $n \to \infty$, it follows that $c_n$ is monotonically decreasing to  $\sqrt{2}$.   Therefore maximum value of $c_n$ occurs at $c_0$.

The open loop system is asymptotically stable because all its eigenvalues are in the open left half plane.
Let ${\bf Z}^o\subset {\bf L}_2[0,1]$ be the closure of the span of $\phi_n(x) $ for $i=0,1,2,\ldots$.
The superscript $^o$ denotes that this is the closure of the  domain of the open loop operator $A$.
This operator is densely defined  on ${\bf Z}^o$ and  generates a strongly continuous semigroup.  If $z^0(x)\in {\bf Z^o}$
then  
\bean
z(x,t)&=& T^{o}(t) z^0(x) \\
&=& \sum_{i=0}^\infty e^{\lambda_n t}  \phi_n(x)\ \int_0^1 \phi_n(x_1)z^0(x_1) \ dx_1 
\eean 
Again the superscript on $T^o$ denotes that this is  the open loop semigroup.

We seek a feedback control law of the form
\bean
u(t)&=& \int_0^1 K(x) z(x,t) \ dx
\eean
to speed up the stabilization.  To find $K(x)$ we solve
a linear quadratic regulator  (LQR). Minimize by choice of $u(t)$ the quantity
\bea \label{crit}
&&\int_0^\infty \iint_{\cal S} Q(x_1,x_2)z(x_1,t)z(x_2,t)\ dA\ dt\\
&& +\int_0^\infty R u^2(t) \nonumber
\eea
subject to (\ref{old}, \ref{ic}, \ref{BC0}, \ref{BC1})
where ${\cal S}$ is the unit square $[0,1]\times[0,1]$ and $dA=dx_1dx_2$.

We require that $R>0$ and $Q(x_1,x_2) $ is a symmetric function, $Q(x_1,x_2) =Q(x_2,x_1), $  satisfying
\bean
0&\le& \iint_{\cal S} Q(x_1,x_2)\theta(x_1)\theta(x_2)\ dA
\eean
for any function $\theta(x)$.  We allow  $Q(x_1,x_2) $ to be a generalized function.   For example 
if $Q(x_1,x_2) =Q(x_1) \delta(x_1-x_2)$, where $ Q(x_1)\ge 0 $ for each $x_1$ and $\delta(x_1-x_2)$ is the  Dirac $\delta$ function then 
\bean
\iint_{\cal S} Q(x_1,x_2)z(x_1,t)z(x_2,t)\ dA= \int_0^1 Q(x) z^2(x,t)\ dx
\eean

Let $P(x_1,x_2)$ be any continuous, symmetric function of $(x_1,x_2)$.  Assume that the
 control trajectory $u(t)$ results in $z(x,t)\to 0$ as $t\to \infty$. We know that 
 such control trajectories exist because the open loop rod, $u(t)=0$, is asymptotically stable.
  By the Fundamental
Theorem of Calculus
\bean
0&=&  \iint_{\cal S} P(x_1,x_2)z^0(x_1)z^0(x_2)\ dA\\ 
&+&\int_0^\infty\iint_{\cal S} {d\over dt}\left(P(x_1,x_2)z(x_1,t)z(x_2,t)\right)\ dA\ dt
\eean
so
\bea \label{P1}
&&0=  \iint_{\cal S}  P(x_1,x_2)z^0(x_1)z^0(x_2)\ dA\\ 
&&+\int_0^\infty \iint_{\cal S}  P(x_1,x_2)\frac{\partial^2 z}{\partial x_1^2}(x_1,t)z(x_2,t) \nonumber\\
&&+ 
P(x_1,x_2)z(x_1,t)\frac{\partial^2 z}{\partial x_2^2}(x_2,t)
\ dA\ dt \nonumber
\eea

We assume that $P(x_1,x_2)$ satisfies  Neumann boundary conditions at $x_i=0$ 
\bea  \label{BC01}
\frac{\partial P}{\partial x_1}(0,x_2)&=&0\\
\frac{\partial P}{\partial x_2}(x_1,0)&=&0  \label{BC02}
\eea
and 
Robin boundary conditions at $x_i=1$
\bea \label{BC11}
\frac{\partial P}{\partial x_1}(1,x_2)&=&-\beta  P(1,x_2)\\
\frac{\partial P}{\partial x_2}(x_1,1)&=&-\beta  P(x_1,1)
\label{BC12}
\eea
Because of the symmetry of $P(x_1,x_2)$, (\ref{BC01}) is equivalent to (\ref{BC02})
and (\ref{BC11}) is equivalent to (\ref{BC12}).

When we integrate (\ref{P1}) by parts twice we get the equation
\bea \label{P2}
0&= & \iint_{\cal S} P(x_1,x_2)z^0(x_1)z^0(x_2)\ dA\\ 
&+&\int_0^\infty  \iint_{\cal S} \nabla^2P(x_1,x_2)z(x_1,t)z(x_2,t)\ dA\  dt\nonumber\\
&+& \int_0^\infty \int_0^1   P(1,x_2)\beta u(t)z(x_2,t)\ dx_2\ dt\nonumber\\
&+& \int_0^\infty \int_0^1  P(x_1,1)z(x_1,t)\beta u(t)\ dx_1\ dt\nonumber\
\eea
where $\nabla^2$ is the two dimensional Laplacian.

We add the right side of (\ref{P2}) to the criterion (\ref{crit}) to be minimized to get the equivalent criterion
\bea \label{P3}
&&  \iint_{\cal S}  P(x_1,x_2)z^0(x_1)z^0(x_2)\ dA\\ 
&+&\int_0^\infty \iint_{\cal S}  \left(\nabla^2P(x_1,x_2)+Q(x_1,x_2)\right)\nonumber\\
&\times&z(x_1,t)z(x_2,t)\ dA\  dt\nonumber\\
&+& \int_0^\infty \int_0^1  P(1,x_2)\beta u(t) z(x_2,t)\ dx_2\ dt\nonumber\\
&+& \int_0^\infty \int_0^1  P(x_1,1) z(x_1,t)\beta u(t)\ dx_1\ dt\nonumber\\
&+& \int_0^\infty  Ru^2(t)
\ dt \nonumber
\eea

Suppose  this criterion can be written as a double integral with respect to the initial state and a time integral of a perfect square involving the control.  If the latter can be made zero by the proper choice of the control then the optimal cost is the double integral with respect to the initial state.

So we would like to chose $K(x)$ so that the time integrand in (\ref{P3}) is  a perfect square.  In other words we want (\ref{P3}) to be of the form 
\bea \label{P4}
&&0=  \iint_{\cal S}  P(x_1,x_2)z^0(x_1)z^0(x_2)\ dA\\ 
&&+ \int_0^\infty \iint_{\cal S}  R\left(u(t)-K(x_1)z(x_1,t)\right)\nonumber\\
&& \times \left(u(t)-K(x_2)z(x_2,t)\right)\ dA\ dt \nonumber 
\eea
Clearly the terms quadratic in $u(t) $ match so we equate the terms involving $u(t)$ and $z(x_1,t)$,
\bean
&& \int_0^\infty \int_0^1 P(x_1,1) z(x_1,t)\beta u(t) \ dx_1dt\\
&&= -\int_0^\infty\int_0^1 K(x_1)Rz(x_1,t) u(t) \ dx_1 dt
\eean
To make these equal we set 
\bea \label{K1}
K(x_1)=-   \beta R^{-1}  P(x_1,1)
\eea
  By the symmetry of $P(x_1,x_2)$,
$K(x_2)=-    \beta R^{-1} P(1,x_2)$.

Finally we equate the terms involving $z(x_1,t)$ and $z(x_2,t)$,
\bean
&& \iint_{\cal S}  \left(\nabla^2P(x_1,x_2)+Q(x_1,x_2)\right)\\
&&\times z(x_1,t)z(x_2,t)\ dA \ dt\\
&&=R\int_0^\infty  \iint_{\cal S} K(x_1)K(x_2)z(x_1,t)z(x_2,t)\ dA\ dt
\eean

This yields what we call a Riccati PDE 
\bea \label{RPDE}
&&\nabla^2P(x_1,x_2)+Q(x_1,x_2)=RK(x_1) K(x_2)\nonumber\\
&&=   \beta^2 R^{-1}    P(x_1,1)P(1,x_2)
\eea
Since we only assumed that $P(x_1,x_2)$ is continuous and $Q(x_1,x_2)$ might equal $Q(x_1)\delta(x_1-x_2)$, this is to be interpreted in the weak sense, if $\theta(x)$ is $C^2$ on $ 0\le x\le 1$ then 
\bean
&&\iint_{\cal S}  \left(\nabla^2P(x_1,x_2)+Q(x_1,x_2)\right)\theta(x_1)\theta(x_2)\ dA\\
&&= \iint_{\cal S}    \beta^2 R^{-1}     P(x_1,1) P(1,x_2)\theta(x_1)\theta(x_2)\ dA
\eean
where the evaluations are done using integration by parts.

The boundary conditions (\ref{BC01}, \ref{BC02} \ref{BC11}, \ref{BC12}) are also to be interpreted in the weak sense
\bean
0&=&\int_0^1  \frac{\partial P}{\partial x_1} (0,x_2) \theta(x_2)\ dx_2\\
0&=&\int_0^1  \frac{\partial P}{\partial x_2} (x_1,0) \theta(x_1)\ dx_1\\
0&=&\int_0^1 \left(\frac{\partial P}{\partial x_1} (1,x_2)+ \beta R^{-1} P(1,x_2)\right) \theta(x_2)\ dx_2\\
0&=&\int_0^1  \left(\frac{\partial P}{\partial x_2} (x_1,1) +  \beta R^{-1}P(x_1,1)\right)\theta(x_1)\ dx_1
\eean

If we can solve the Riccati PDE  subject to these boundary conditions then clearly the optimal cost starting 
from $z^0(x)$ is a quadratic functional of the initial state 
\bean
\iint_{\cal S}  P(x_1,x_2)z^0(x_1)z^0(x_2) \ dA
\eean
and the optimal  feedback is a linear  functional of the current state
\bea 
u(t)&=& \int_0^1K(x) z(x,t) \ dx  \nonumber\\
&=&-   \beta R^{-1}  \int_0^1P(x,1) z(x,t) \ dx \label{olf}
\eea

We assume that the solution to the Riccati PDE has an expansion in the open loop eigenfunctions (\ref{ef}),
\bea \label{Riccexp}
P(x_1,x_2)&=& \sum_{n_1,n_2=0}^\infty P_{n_1,n_2}
\phi_{n_1}(x_1)\phi_{n_2}(x_2)
\eea
In abuse of notation we use the symbol $P$ to denote both a function $P(x_1,x_2)$ and a coefficient $P_{n_1,n_2}$.  The proper meaning 
should be clear from context.
Clearly any such expansion satisfies the boundary conditions (\ref{BC01}, \ref{BC02}, \ref{BC11}, \ref{BC12}).
Then
\bea \label{RiccexpK}
K(x)&=& - \beta R^{-1} \sum_{n_1,n_2=0}^\infty P_{n_1,n_2}
\phi_{n_1}(1)\phi_{n_2}(x)
\eea
 Because we seek a symmetric weak solution without
loss of generality   $ P_{n_1,n_2}= P_{n_2,n_1}$.

We also assume that $Q(x_1,x_2)$ has a similar expansion,
\bea \label{Qexp}
Q(x_1,x_2)&=& \sum_{i_1,i_2=0}^\infty Q_{n_1,n_2}\phi_{n_1}(x)\phi_{n_2}(x)
\eea
We plug these into (\ref{RPDE}) and get an algebraic Riccati equation for the infinite dimensional 
matrix $ P_{n_1,n_2}$,
\bea
&&(\lambda_{n_1}+\lambda_{n_2}) P_{n_1,n_2} + Q_{n_1,n_2}\label{mre} \\
&&=   \beta^2 R^{-1}   \sum_{m_1,m_2=0}^\infty  P_{n_1,m_1}  P_{n_2,m_2} \phi_{m_1}(1)\phi_{m_2}(1)\nonumber
\eea

To simplify the notation henceforth we assume 
  \bea 
  \beta &=&1 \label{asmp}\\
   Q(x_1,x_2)&=&\delta(x_1-x_2)\nonumber\\
  R&=&1\nonumber
  \eea
  where  $\delta(x_1-x_2)$ is the Dirac $\delta$. 
  Then by Parseval's equality
  \bean
  Q_{n_1,n_2}&=& \delta_{n_1,n_2} \nonumber
\eean
where $ \delta_{n_1,n_2}$ is the Kronecker $\delta$.

  Let $\Pi$ denote the infinite dimensional matrix $[ P_{n_1,n_2}]$.
  Then (\ref{mre}) is the algebraic Riccati equation 
  \bean
  F'\Pi+\Pi F +Q&=& \Pi G R^{-1}G' \Pi
  \eean
of the infinite dimensional linear quadratic  system where
  \bean
  F&=&\bmt \lambda_1&0&0& \ldots\\
  0&\lambda_2&0&\ldots\\
  0&0& \ddots &
  \emt\\
  G&=& \bmt \phi_1(1)\\ \phi_2(1) \\  \vdots \emt\\
  Q&=&\bmt 1&0&0& \ldots\\
  0&1&0&\ldots\\
  && \ddots &
  \emt\\
  R&=&1
  \eean

 We denote our first guess at a solution to (\ref{mre}) by $ P^{(1)}_{n_1,n_2}$   and we take it to be  diagonal,
   \bea
   \label{diag}
   P^{(1)}_{n_1,n_2}=\delta_{n_1,n_2} P^{(1)}_{n_1,n_2}
  \eea 
Then we get a sequence of quadratic equations for $ P^{(1)}_{n,n}$,
 \bea \label{dr}
0&=&   \phi^2_{n}(1)  \left( P^{(1)}_{n,n}\right)^2- 2\lambda_{n,n} P^{(1)}_{n,n}-1
\eea
Clearly we wish to take the positive root so we assume
\bea \label{Pi2}
 P^{(1)}_{n,n}={ \lambda_{n}+\sqrt{\lambda_{n}^2+  \phi_{n}^2(1)}\over   \phi_{n}^2(1)}
\eea

But  we need to check if the off-diagonal terms  satisfy  (\ref{mre}).    
   If (\ref{diag}) holds then when $n_1\ne n_2$  (\ref{mre}) becomes
\bean
0&=& \sum_{m_1,m_2=0}^\infty  P^{(1)}_{n_1,m_1}  P^{(1)}_{n_2,m_2} \phi_{m_1}(1)\phi_{m_2}(1)\\
&=& P^{(1)}_{n_1,n_1}  P^{(1)}_{n_2,n_2} \phi_{n_1}(1)\phi_{n_2}(1)\ \neq\ 0
\eean
So we conclude that the solution to  (\ref{mre}) is not diagonal.   

Given $ P^{(k)}_{n_1,n_2}$ we define $ P^{(k+1)}_{n_1,n_2}$ by the equation
 \bea
&&(\lambda_{n_1}+\lambda_{n_2}) P^{(k+1)}_{n_1,n_2} + Q_{n_1,n_2}\label{mrek} \\
&&=      \sum_{m_1,m_2=0}^\infty  P^{(k)}_{n_1,m_1}  P^{(k)}_{n_2,m_2} \phi_{m_1}(1)\phi_{m_2}(1)\nonumber
\eea
Then we define
\bean
P^{(k+1)}(x_1,x_2)&=& \sum_{n_1,n_2=0}^\infty P^{(k+1)}_{n_1,n_2}
\phi_{n_1}(x_1)\phi_{n_2}(x_2)\\
K^{(k+1)}(x)&=& \sum_{n_1,n_2=0}^\infty P^{(k+1)}_{n_1,n_2}
\phi_{n_1}(1)\phi_{n_2}(x)
\eean



Although it may not be obvious this is variation of the familiar policy iteration scheme for solving an optimal control problem.  Our $k^{th} $ approximation of the optimal cost of starting at $z^0(x)$ is 
\bean
\iint_{\cal S}  P^{(k)}(x_1,x_2)z^0(x_1)z^0(x_2) \ dA
\eean
where
\bean
 P^{(k)}(x_1,x_2)&=& \sum_{n_1,n_2=0}^\infty P^{(k)}_{n_1,n_2} \phi_{n_1}(x_1)\phi_{n_2}(x_2)
 \eean
 Given this approximation then we plug $P^{(k)}(x_1,x_2)$ into (\ref{P3}) to get
\bea \label{P3k}
&&  \iint_{\cal S}  P^{(k)}(x_1,x_2)z^0(x_1)z^0(x_2)\ dA\\ 
&+&\int_0^\infty \iint_{\cal S}  \left(\nabla^2P^{(k)}(x_1,x_2)+\delta(x_1-x_2)\right)\nonumber\\
&\times&z(x_1,t)z(x_2,t)\ dA\  dt\nonumber\\
&+& \int_0^\infty \int_0^1  P^{(k)}(1,x_2) u(t) z(x_2,t)\ dx_2\ dt\nonumber\\
&+& \int_0^\infty \int_0^1  P^{(k)}(x_1,1) z(x_1,t) u(t)\ dx_1\ dt\nonumber\\
&+& \int_0^\infty  u^2(t)
\ dt \nonumber
\eea
To find  the $k^{th} $ approximation of the optimal control we minimize this expression with respect to $u$ and
obtain
\bean
u^{(k)}(t)&=&-\int_0^1  P^{(k)}(x_1,1) z(x_1,t)\ dx_1
\eean 
so the $k^{th} $ approximation of the optimal gain is 
\bean
K^{(k)}(x_1)&=& -P^{(k)}(x_1,1)
\eean
Let $z^{k}(x,t)$ be the solution of (\ref{old}, \ref{ic}, \ref{BC0}, \ref{BC1}) when $u(t)=u^{(k)}(t)$.
The $(k+1)^{th} $ approximation of the optimal cost is
\bean
&&  \iint_{\cal S}  P^{(k)}(x_1,x_2)z^0(x_1)z^0(x_2)\ dA\\
&&=\int_0^\infty \int_0^1 \left(z^{k}(x,t)\right)^2+\left(u^{(k)}(t)\right)^2 \ dx\ dt
\eean

Then for any $z^0(x)$ the sequence of scalars
\bean
\iint_{\cal S}  P^{(k)}(x_1,x_2)z^0(x_1)z^0(x_2)\ dA
\eean
is nonincreasing in $k$ and bounded below by zero hence it is convergent.  It follows that the Fourier
coefficients $P^{(k)}_{n_1,n_2}$ are also convergent.

We compute an approximation
 to the  upper left block of $ P^{(50)}_{n_1,n_2} $ by
 truncating the iteration (\ref{mrek}) to
 \bean
&&(\lambda_{n_1}+\lambda_{n_2}) P^{(k+1)}_{n_1,n_2} + Q_{n_1,n_2} \\
&&=      \sum_{m_1,m_2=0}^{10}  P^{(k)}_{n_1,m_1}  P^{(k)}_{n_2,m_2} \phi_{m_1}(1)\phi_{m_2}(1)\nonumber
\eean
Then after $50$ iterations the upper left block of $ P^{(50)}_{n_1,n_2} $ is approximately
\bean
\bmt
 0.5757  & -0.0018  & -0.0002  & -0.0000    \\
   -0.0018  &  0.0425  & -0.0000  & -0.0000   \\
   -0.0002 &  -0.0000  &  0.0121 &  -0.0000   \\
   -0.0000  & -0.0000  & -0.0000  &  0.0055  
   \emt
   \eean
Notice how strongly diagonally dominant this is and how the diagonal elements are decreasing quite fast.\\

\section{Closed Loop Eigenvalues and Eigenvectors}
We continue to assume that (\ref{asmp}) holds.
  The boundary feedback does not appear in the closed loop dynamics, it is still 
  \bean
  \frac{\partial z}{\partial t}(x,t)&=&    \frac{\partial^2 z}{\partial x^2}(x,t)
  \eean
but the boundary conditions are changed by the feedback.  The boundary
condition at $x=0$ is still the Neumann boundary
condition  (\ref{BC0}) but the Robin boundary
condition (\ref{BC1})  at $x=1$ is replaced by
\bea \label{BCC1}
&&\frac{\partial z}{\partial x}(1,t)= \int_0^1K(x_1)z(x_1,t)\ dx_1-z(1,t)\\
&& =\sum_{n_1,n_2=0}^\infty P_{n_1,n_2}\phi_{n_2}(1)\int_0^1\phi_{n_1}(x)z(x,t)\ dx-z(1,t) \nonumber
\eea
Note that this nonstandard boundary condition (\ref{BCC1}) is linear in $z(x,t)$.

Because of the Neumann BC at $x=0$ we know that the unit normal closed loop eigenvectors are of the form
\bea  \label{psi}
\psi(x)&=&c(\rho)\cos \rho x
\eea
for some $\rho>0$
where the normalizing constant  is given by 
\bea  
c(\rho)&=&\mbox{sign}(\cos \rho)\sqrt{{4\rho\over 2\rho +\sin 2\rho}} 
\eea

The $ \rho $ are chosen so that $g(\rho)=0$ where
\bea \label{rheq}
&&g(\rho)=\rho \sin \rho-  \cos \rho \\
&& - \sum_{n_1,n_2=0}^\infty P_{n_1,n_2}\phi_{n_2}(1)\int_0^1\phi_{n_1}(x_1)\cos \rho x_1\ dx_1 \nonumber
\eea
 The one dimensional Laplacian is not self adjoint under these boundary conditions (\ref{BC0}, \ref{BCC1}) 
so there is no reason to expect that the closed loop eigenfunctions are orthogonal.\\

   If $n>10$ then  $ P_{n,n}\le {1\over 2 |\lambda_{10}|}\le  {1\over 200 \pi^2}\approx  0.0005$, so
to approximate  the first few $\rho_n$ we truncate (\ref{rheq}) to 
\bea \label{rheq2}
&&\rho \sin \rho-\cos \rho\\
&& =  \sum_{n_1.n_2=0}^{10} P_{n_1,n_2}\phi_{n_2}(1)\int_0^1\phi_{n_1}(x_1)\psi(x_1)\ dx_1 \nonumber
\eea
We solve (\ref{rheq2}) by Newton's method starting at $\rho =1.05\nu_n$
for $n=0,1,2,3,4$.  The result is  
  $\rho_0\approx 0.9982$, $\rho_1\approx  3.4381$, $\rho_2\approx 6.4391$,  $\rho_3\approx 9.5299$  and $\rho_4\approx 12.6455$.
  
The first five closed loop eigenvalues are approximately
$\mu_0\approx -\rho_0^2=-0.9964$, $\mu_1\approx-\rho_1^2=-11.8202$, $\mu_2\approx -\rho_2^2= -41.4618$,  $\mu_3\approx-\rho_3^2= -90.8190$ and $\mu_4\approx-\rho_4^2= -159.9095$.  
Recall 
the first five open loop eigenvalues are approximately $\lambda_0\approx-0.7401$, $\lambda_1\approx-11.7347$, $\lambda_2\approx-41.4620$, $\lambda_3\approx -90.8192$  and $\lambda_4\approx -159.9095$.  Notice how close $\mu_n$ and $\lambda_n$ are if $n>0$.
   The boundary feedback has  
a significant   effect on the least stable open loop  eigenvalue but less so on the rest of the open loop eigenvalues because they are already so stable it would cost to much control energy to significantly increase their stability.  

\section{ Boundary Control of a Nonlinear  Reaction Diffusion Equation}
To the above system we add a destabilizing nonlinear term to obtain the boundary controlled reaction diffusion system
\bea
\frac{\partial z}{\partial t}(x,t)&=&\frac{\partial^2 z}{\partial x^2}(x,t)+\alpha z^2(x,t) \label{nlrdyn}\\
z(x,0)&=& z^0(x) \label{icc}\\ 
\frac{\partial z}{\partial x}(0,t)&=&0 \label{bcc0}\\
\frac{\partial z}{\partial x}(1,t)&=&  (u(t)-z(1,t))\label{bcc1}
\eea
for some positive constant $\alpha$.   Vazquez and Krstic \cite{VK08} used backstepping to stabilize a similar system  with $\alpha=1$.
They assumed a different boundary condition at $x=1$, namely direct control of the heat flux,
\bean
\frac{\partial z}{\partial x}(1,t)&=&u(t)
\eean

To find  a feedback to stabilize this system
we consider the nonlinear quadratic optimal control of minimizing (\ref{crit}) subject to (\ref{nlrdyn}, \ref{icc}, 
\ref{bcc0}, \ref{bcc1}).

It is well known \cite{FZ00} that this system cannot be globally stabilized but we are only interested  in local stabilization
around $z=0$.   The reason is that this is a mathematical model of a physical system and the model is not globally valid, there
is an absolute zero temperature that the rod cannot go below and at a sufficiently high temperature the rod will melt.
So global stabilization is of mathematical but not of physical interest.

Let $P(x_1,x_2)$ be the solution of the Riccati PDE (\ref{RPDE}) and $K(x)$ be the gain of the optimal linear feedback (\ref{K1}).
In abuse of notation let $P(x_1,x_2,x_3) $ be a symmetric function of three variables. We distinguish between $P(x_1,x_2)$ and $P(x_1,x_2,x_3)$ by the number of arguments. Symmetric  means that the value of the function is invariant
under any permutation of the three variables. 

We further assume that  $P(x_1,x_2,x_3) $
weakly satisfies  Neumann boundary conditions at $x_1=0$ 
\bea \label{BC03}
\frac{\partial P}{\partial x_1}(0,x_2,x_3)&=&0
\eea
and Robin boundary conditions at $x_1=1$
\bea \label{BC13}
\frac{\partial P}{\partial x_1}(1,x_2,x_3)&=&-  P(1,x_2,x_3)  
\eea
By symmetry, similar boundary conditions hold at $x_2=0,1$ and $x_3=0,1$.
  
    Assume also that the optimal feedback takes the form
\bea \label{k12}
u(t)&=& \int_0^1 K(x_1) z(x_1,t)\ dx_1 
\\&&+\iint_{\cal S}  K(x_1,x_2) z(x_1,t)z(x_2,t)\ dA \\
&&+O(\left\|z(x,t)\right\|)^3 \nonumber
\eea
where
 \bean
\left\|z(x,t)\right\|^2&=& \int_0^1|z(x,t)|^2\ dx
 \eean 
   Again we distinguish between the
optimal linear feedback gain $K(x_1)$ and the optimal quadratic feedback gain $K(x_1,x_2)$ by the number of arguments

Again by the Fundamental Theorem of Calculus if the control trajectory $u(t)$  takes $z(x,t)\to 0$
as $t \to \infty$ then
\bea \label{P333}
&&0=  \iint_{\cal S}P(x_1,x_2)z^0(x_1) z^0(x_2)\ dA\\ 
&&+ \iiint_{\cal C} P(x_1,x_2,x_3)z^0(x_1)z^0(x_2)z^0(x_3)\ dV\nonumber \\
&&+\int_0^\infty \iint_{\cal S}{d\over dt}\left(P(x_1,x_2)z(x_1,t)z(x_2,t)\right)\ dA\ dt\nonumber\\
&&+\int_0^\infty \iiint_{\cal C} {d\over dt}\left(P(x_1,x_2,x_3)z(x_1,t)z(x_2,t)z(x_3,t)\right)\nonumber\\
&& \times  dV \ dt\nonumber+O(\left\|z(x,t)\right\|)^4
\eea
where ${\cal C}$ denotes the unit cube $[0,1]\times[0,1]\times[0,1]$ and $dV$ is the volume element
$dV=dx_1dx_2dx_3$.

Because $P(x_1,x_2)$ is the solution of the Riccati PDE, the terms quadatic in $z$ in the time integral drop out.  But we pick up  cubic terms
from the boundary when we integrate (\ref{P333}) by parts twice.  We use the symmetry  of $P(x_1,x_2) $ and $P(x_1,x_2,x_3) $
to condense them. Then we obtain 
\bea \label{P33}
&&0=   \iint_{\cal S}P(x_1,x_2)z^0(x_1) z^0(x_2)\ dA \\  &&
+\int_0^\infty \iiint_{\cal C} P(x_1,x_2,x_3)z^0(x_1)z^0(x_2)z^0(x_3)\ dV \nonumber \\
&&+2\alpha \int_0^\infty\iint_{\cal S} P(x_1,x_2)  z^2(x_1,t)z(x_2,t) \ dA\  dt\nonumber\\
&&+3\int_0^\infty\iiint_{\cal C} P(x_1,x_2,x_3)\frac{\partial^2 z}{\partial x_1^2}(x_1,t)\nonumber\\
&&\times z(x_2,t)z(x_3,t)  dV\ dt \nonumber\\
&&+3 \int_0^\infty \iiint_{\cal C} P(1,x_2,x_3) K(x_1)\nonumber\\
&&\times z(x_1,t)z(x_2,t) z(x_3,t)\ dV\ dt\nonumber \\
&&+2  \int_0^\infty \iiint_{\cal C} P(x_1,1) K(x_2,x_3)\nonumber\\
&& \times z(x_1,t)z(x_2,t) z(x_3,t)\ dV\ dt +O(\left|z(x,t)\right|)^4 \nonumber 
\eea

If we integrate (\ref{P33}) by parts twice 
we get the equation
\bea \label{P34}
&&0=  \iint_{\cal S} P(x_1,x_2)z^0(x_1) z^0(x_2)\ dA\\ 
&&+ \iiint_{\cal C} P(x_1,x_2,x_3) z^0(x_1)z^0(x_2)z^0(x_3)\ dV\nonumber\\ 
&&+\int_0^\infty \iiint_{\cal C} \nabla^2P(x_1,x_2,x_3)\nonumber \\
&&\times z(x_1,t)z(x_2,t)z(x_3,t)\  dV dt\nonumber\\
&&+2\alpha \int_0^\infty \iint_{\cal S} P(x_1,x_2) z^2(x_1,t)z(x_2,t)\ dA\ dt\nonumber
\\
&&+3 \int_0^\infty \iiint_{\cal C} P(1,x_2,x_3) K(x_1) \nonumber
 \\&&\times z(x_1,t)z(x_2,t)z(x_3,t)\  dV dt\nonumber\\
&&+2 \int_0^\infty \iiint_{\cal C}P(x_1,1) K(x_2,x_3)\nonumber\\
&& \times z(x_1,t)z(x_2,t) z(x_3,t)\ dVdt +O(\left\|z(x,t)\right\|)^4 \nonumber 
\eea
where $\nabla^2$ is now the three dimensional Laplacian.  This equation (\ref{P34})
is not symmetric in $x_1,x_2,x_3$ but we are looking for a symmetric  weak  solution.  If 
$P(x_1,x_2,x_3)$ is a weak solution that is not symmetric then we can get a symmetric weak
solution by averaging over all permutations of $x_1,x_2,x_3$.

We have cancelled the quadratic terms in the criterion by our choice $P(x_1,x_2)$ and $K(x)$ but the quadratic term in the feedback generates a cubic term in the criterion of the form
\bean
 &&2\int_0^\infty \int_0^1 \int_0^1\int_0^ 1 K(x_1) K(x_2,x_3) \\
 && \times z(x_1,t)z(x_2,t) z(x_3,t)\ dV\ dt 
   \eean
  Because of (\ref{K1}) this equals
  \bean
 && -2 \int_0^\infty \iiint_{\cal C}P(x_1,1) K(x_2,x_3)\nonumber\\
&& \times z(x_1,t)z(x_2,t) z(x_3,t)\ dV\ dt 
\eean
and so this cancels out the last term in (\ref{P34}).  

So the new criterion to be minimized is 
\bea \label{P35}
&&0=  \iint_{\cal S} P(x_1,x_2)z^0(x_1) z^0(x_2)\ dA\\ 
&&+ \iiint_{\cal C} P(x_1,x_2,x_3) z^0(x_1)z^0(x_2)z^0(x_3)\ dV \nonumber
\\&&+\int_0^\infty \iiint_{\cal C}  \nabla^2P(x_1,x_2,x_3)\nonumber\\
&& \times z^0(x_1)z^0(x_2)z^0(x_3)\ dV \ dt\nonumber\\
&&+3  \int_0^\infty \iiint_{\cal C} P(x_1,x_2,1) K(x_3) \nonumber\\
&&\times z(x_1)z(x_2)z(x_3)
\ dV\ dt \nonumber\\ 
&&+2\alpha \int_0^\infty \iint_{\cal S} P(x_1,x_2) z^2(x_1,t)z(x_2,t)\ dA\ dt\nonumber
\eea
Notice the quadratic optimal gain $K(x_1,x_2)$ does not appear in this equation so we can solve  for
cubic kernel $P(x_1,x_2,x_3)$ of the optimal cost independently of $K(x_1,x_2)$.

We assume that $P(x_1,x_2,x_3) $ is a symmetric weak solution to the symmetric  linear elliptic PDE
\bea \label{3pde}
&&0=\nabla^2P(x_1,x_2,x_3)
\\&&
+  P(1,x_2,x_3) K(x_1)\nonumber\\
&&+P(x_1,1,x_3) K(x_2)\nonumber\\
&&+P(x_1,x_2,1) K(x_3)\nonumber
\nonumber \\ &&+{2\alpha\over 3} \left(P(x_1,x_2)\delta(x_1-x_3)+P(x_2,x_3)\delta(x_2-x_1)\right. \nonumber\\
&& \left. +P(x_3,x_1)\delta(x_3-x_2)\right)\nonumber
\eea
subject to the weak boundary conditions (\ref{BC03}, \ref{BC13}).
Then the optimal cost is
\bea \label{oc3}
&& \iint_{\cal S}P(x_1,x_2)z^0(x_1) z^0(x_2)\ dA\\ 
&&+ \iiint_{\cal C}P(x_1,x_2,x_3)z^0(x_1)z^0(x_2)z^0(x_3)\ dV\nonumber \\&& +O(\left\|z^0(x)\right\|^4)\nonumber
\eea
where
\bean
\left\|z^0(x)\right\|^2=\int_0^1\left(z^0(x)\right)^2\ dx
\eean
 
To find the quadratic  optimal feedback gain we
 start by noting that the  optimal cost starting from $z^0(x)$ is (\ref{oc3}).

 Let $u^*(t) $ be the optimal control trajectory starting from $z^0(x)$ then by arguments similar to the above we obtain
\bean
&&\int_0^\infty\int_0^1 z^2(x,t)\ dx +(u^*(t))^2\ dx\  dt\\
&&+\int_0^\infty \int_{[0,1]^2}\left( \nabla^2P(x_1,x_2)+2\alpha_1\right)\\
&& \times z(x_1,t)z(x_2,t) \ dA\ dt\\
&&+\int_0^\infty \int_{[0,1]^3} \nabla^2P(x_1,x_2,x_3)\\
&& \times z(x_1,t)z(x_2,t)z(x_3,t) \ dV\ dt\\
&&+2\alpha_2 \int_0^\infty \int_{[0,1]^2} \nabla^2P(x_1,x_2)z^2(x_1,t)z(x_2,t)\\
&&+2\int_0^\infty\int_0^1P(x_1,1)z(x_1,t)u^*(t)\ dx_1\ dt\\
&&+3\int_0^\infty \int_{[0,1]^2} P(x_1,x_2,1)\\
&& \times z(x_1,t)z(x_2,t)u^*(t) \ dA\ dt\\
&&+O(\left\|z^0(x)\right\|^4)
\eean

 We replace $u^*(t)$ by $u^*(t)+\epsilon \delta u(t)$ in this expression and  differentiate
with respect to $\epsilon$ at $\epsilon=0$ to obtain the expression 
 \bean
&&2\int_0^\infty u^*(t)\delta u(t)\ dt\\
&&+2\int_0^\infty\int_0^1P(x_1,1)z(x_1,t)\delta u(t)\ dx_1\ dt\\
&&+3\int_0^\infty \int_{[0,1]^2} P(x_1,x_2,1)z(x_1,t)z(x_2,t)\delta u(t) \ dA\ dt\\
&&+O(\left\|z^0(x)\right\|^4)
\eean
Since $u^*(t) $ is the optimal control trajectory this quantity must be zero for 
any $\delta u(t)$ so the coefficient of $\delta u(t)$ must vanish for each
$t$. This leads to the expressions for the linear and quadratic parts of the
optimal feedback,
\bean 
u^*(t)&=&\int_0^1K(x_1)z(x_1,t)\ dx_1\\
&&+\int_{[0,1]^2} K(x_1,x_2)z(x_1,t)z(x_2,t)\ dA\\ 
&&+O(\left\|z^0(x)\right\|^3)
\eean
where
\bea \label{K11} 
K(x_1)&=&-P(x_1,1)\\
K(x_1,x_2)&=&-{3\over 2}P(x_1,x_2,1)
\label{K2} 
\eea

Again we expand  $P(x_1,x_2,x_3)$ in terms of the orthonormal eigenfunctions of the unforced system
\bean 
P(x_1,x_2,x_3)= \sum_{n_1,n_2,n_3=1}^\infty P_{n_1,n_2,n_3}\phi_{n_1}(x_1)\phi_{n_2}(x_2)\phi_{n_3}(x_3)
\eean
where $P_{n_1,n_2,n_3}$ is a three tensor that is symmetric in its three indices.  The symmetric linear elliptic PDE
(\ref{3pde}) is then equivalent to
\bea
&&0=\left(\lambda_{n_1}+\lambda_{n_2}+\lambda_{n_3}\right)P_{n_1,n_2,n_3}\\
&&-\sum_{m_1,m_2=1}^\infty P_{m_1,n_2,n_3}P_{m_1,m_2}\phi_{m_2}(1) \nonumber\\
&&-\sum_{m_1,m_2=1}^\infty P_{n_1,m_1,n_3}P_{m_1,m_2}\phi_{m_2}(1) \nonumber\\
&&-\sum_{m_1,m_2=1}^\infty P_{n_1,n_2,m_1}P_{m_1,m_2}\phi_{m_2}(1) \nonumber\\
&&+{2\alpha\over 3}\left(P_{n_1,n_2}\delta_{n_1,n_3}+P_{n_2,n_3}\delta_{n_2,n_1}
+P_{n_3,n_1}\delta_{n_3,n_2}\right)\nonumber
\eea

Since $\lambda_n$ is going to $-\infty$ like $-n^2$  and the $P_{n_1,n_2}$ are going to zero like ${1\over n_1^2+n_2^2}$ it follows that the  $P_{n_1,n_2,n_3}$ are going to zero  at least as fast as ${1\over 3n^2}$.  
Furthermore  since $P_{n_1,n_1}$ is much larger than $P_{n_1,n_ 2}$ when $n_1\ne n_2$ it follows that $P_{n_1,n_1,n_1}$ is much larger than the other $P_{n_1,n_2,n_3}$.

From (\ref{K2}) the gain of the quadratic part of the optimal feedback is 
\bean
&&K(x_1,x_2)\\
&& =-{3\over 2}     \sum_{n_1,n_2,n_3=0}^\infty P_{n_1,n_2,n_3} \phi_{n_1}(x_1)\phi_{n_2}(x_2)\phi_{n_3}(1)
\eean
Since $ P_{n_1,n_2,n_3} $ is symmetric in $n_1,n_2,n_3$, it follows that $K(x_1,x_2)$ is symmetric in $x_1,x_2$.

The higher degree terms in the optimal cost and optimal feedback are found in  a similar fashion.

\section{Example} 
We consider the reaction diffusion system (\ref{nlrdyn}, \ref{icc}, \ref{bcc0}, \ref{bcc1})   with $\alpha=1$. We discretize  the system by
 choose a positive integer $n$ and letting $\zeta_k(t) = z(k/n,t)$  for $ k=0,\ldots, n$.  
 The discretization of the differential equation (\ref{nlrdyn})  is
 \bean
 \dot{\zeta}_k&=& {\zeta_{k+1} -2\zeta_k+\zeta_{k-1}\over n^{-2}}+ \zeta_k^2
 \\ &=& n^2(\zeta_{k+1} -2\zeta_k+\zeta_{k-1})+ \zeta_k^2
 \eean
for $k=1,\ldots,n-1$.  

 To handle the boundary conditions, we add ficticious points $\zeta_{-1}$ and $\zeta_{n+1}$.
The Neumann boundary condition (\ref{BC0}) at $x=0$ is discretized
by a centered first difference
\bean
0&=&  \frac{\partial z}{\partial x}(0,t) \approx n{\zeta_1-\zeta_{-1}\over 2} 
\eean
and we solve this for the ficticious point and obtain
$
\zeta_{-1}=\zeta_{1}$.
Then
 \bean
 \dot{\zeta_0}&=& n^2(\zeta_1-2\zeta_0+\zeta_{-1})+ \zeta_0^2\ = \ 2n^2(\zeta_1-\zeta_0)+ \zeta_0^2
 \eean
The controlled boundary condition (\ref{BC1}) at $x=1$ is also discretized
by a centered first difference
\bean
n{\zeta_{n+1}-\zeta_{n-1}\over 2}&\approx&\frac{\partial z}{\partial x}(0,t)  \\
&=&  \left( u- \zeta_n\right)\\
\eean
We solve this for $\zeta_{n+1}$
\bean
\zeta_{n+1}&=&\zeta_{n-1}+  {2  \over n}\left( u- \zeta_n\right)
\eean
  and plug this  into the differential equation for $\zeta_{n}$
  \bean
   \dot{\zeta}_{n}&=& n^2\left(  \zeta_{n+1}-2\zeta_n+\zeta_{n-1}\right)+ \zeta_n^2\\
   &=&2n^2\left(\zeta_{n-1}-{n+\beta\over n}\zeta_n+{\beta\over n}u\right)+ \zeta_n^2
  \eean  
  
The linear part of this  is the $11$ dimensional system
\bean
\dot{\zeta}&=& F\zeta+Gu
\eean
where
\bean
\zeta&=& \bmt \zeta_0&\ldots&\zeta_{n}\emt'
\eean
\bean
F=n^2\bmt -2&2&&&\\
1&-2&1&&&\\
&\ddots&\ddots&\ddots& \\
&& 1&-2& 1\\
&&&2&-{2n+2  \over n}
\emt
\eean
and
\bean
G=\bmt 0&0&\ldots&0&2n\beta  \emt'
\eean

If $n=10$ and $ \beta=1$ the three least stable poles of $F$ are   $-0.7404$, $  -11.6538$ and $-40.1566$.
Recall that $\lambda_0=-0.7402$, $\lambda_1=-11.7349$ and $\lambda_2=-45.3075$ so there is  substantial agreement between the first three open loop poles of the finite and infinite dimensional systems.

We set $Q$ to be $11\times 11$ matrix 
\bean
\bmt 0.5&0&0&\ldots&0\\
0&1&0& \ldots&0\\
&&\ddots&\\
0& \dots && 1 &0\\
0&\ldots&&0&0.5
\emt
\eean
and $R=1$ and  solved the resulting finite dimensional LQR problem.
Its three least stable closed loop poles are $-1.0396$, $-11.7270$ and $ -40.1804$.  
We found previously that the  three least stable closed loop  poles of the infinite dimensional system are $\mu_0\approx -1.0409$, $\mu_1\approx-11.8094$ and $\mu_2\approx-41.4620$, so there is also substantial agreement between the first three closed loop poles of the finite and infinite dimensional systems.

 At the suggestion of Rafael Vazquez  we used the Crank-Nicolson method \cite{CN47} to simulate the nonlinear system.   Crank-Nicolson is an implicit method that uses an average 
of forward and backward Euler steps.   At each time we took a forward Euler step and then used fixed point iteration to correct for backward Euler.  These  converged after $5$ iterations.
The spatial step was $\Delta x=1/10$ and the temporal step was $\Delta t=(\Delta x)^2=1/100$.   We used our Nonlinear Systems Toolbox \cite{NST} to compute the optimal feedback through cubic terms.
It took about $1.5$  seconds  on a MacBook Pro with a Apple M1 Pro processor.

The open loop nonlinear system  converged slowly when the initial condition was $\zeta_i(0)=0.7$ for $i=0,\ldots,10$ but diverged when $\zeta_n(0)=0.8$, see Figure 1.
\begin{figure}
\centering
\includegraphics[width=3.0in]{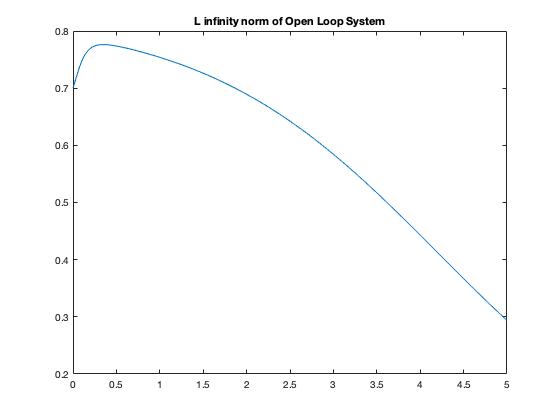}
\caption{L Infinity Norm of the Open Loop Nonlinear System}    
\label{fig:4} 
\end{figure}

The  nonlinear system with optimal linear feedback  converged  when $\zeta_i(0)=1$ for $i=0,\ldots,10$ but diverged when $\zeta_n(0)=1.1$, see Figure 2.
\begin{figure}
\centering
\includegraphics[width=3.0in]{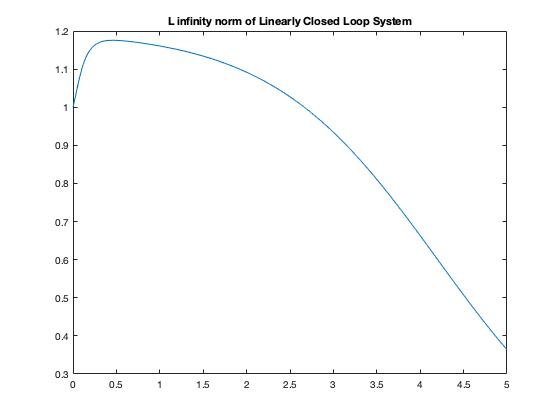}
\caption{L Infinity Norm of the Nonlinear System using Optimal Linear Feedback}    
\label{fig:4} 
\end{figure}

The  nonlinear system  with optimal linear, quadratic and cubic feedback  converged when  $\zeta_i(0)=4$ for $i=0,\ldots,10$ but diverged when $\zeta_n(0)=4.1$, see Figure 3.
\begin{figure}
\centering
\includegraphics[width=3.0in]{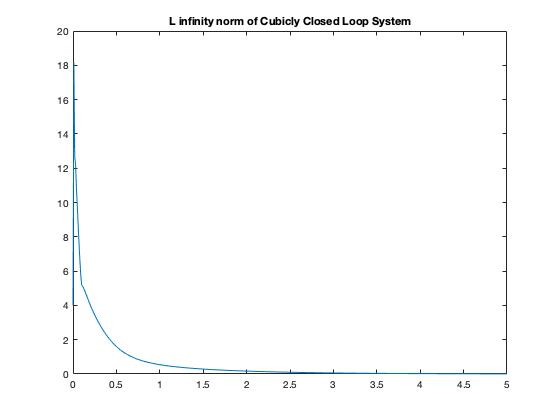}
\caption{L Infinity Norm of the Nonlinear System using Optimal Linear, Quadratic and Cubic Feedback}    
\label{fig:4} 
\end{figure}

\section{Conclusion}
We have solved the LQR problem for the boundary control of  a infinite  dimensional system by extending the finite dimensional technique
of completing the square.  We also optimally locally stabilized a nonlinear reaction diffusion system by extending Al'brekht's method to infinite  dimensions. 
We showed by example  that optimal cubic feedback can lead to a much larger basin of stability than optimal linear feedback.
 We belive that solving an LQR by the completing the square is applicable to other linear infinite dimensional boundary control problems.  We are exploring extending it to 
 the wave equation 
 and the beam equation.

\end{document}